
\font\cc=cmcsc10 scaled 1200 \font\hh=cmcsc10 scaled 1100
  \font\dd=cmr12 
\font\rr=cmr9 \font\ss=cmcsc9 \font\teneufm=eufm10
\font\seveneufm=eufm7 \font\fiveeufm=eufm5
\newfam\eufmfam
\textfont\eufmfam=\teneufm
\scriptfont\eufmfam=\seveneufm
\scriptscriptfont\eufmfam=\fiveeufm
\def\mathfrak#1{{\fam\eufmfam\relax#1}}

\font\tenmsb=msbm10
\font\sevenmsb=msbm7
\font\fivemsb=msbm5
\newfam\msbfam
\textfont\msbfam=\tenmsb
\scriptfont\msbfam=\sevenmsb
\scriptscriptfont\msbfam=\fivemsb
\def\Bbb#1{{\fam\msbfam #1}}

\def \NN {\Bbb N}
\def \CC {\Bbb C}

\def \RR {\Bbb R}

\def\rightheadline{\hfil\rr  On certain large additive functions
\hfil\folio}
\def\leftheadline{\rr\folio\hfil Aleksandar Ivi\'c
 \hfil}
\def\emptyheadline{\hfil}
\headline{\ifnum\pageno=1 \emptyheadline\else
\ifodd\pageno \rightheadline \else \leftheadline\fi\fi}

\def\a{\alpha}
\def\b{\beta} \def\e{\varepsilon}
\def\no{\noindent} \def\d{\,{\rm d}}
\topglue1cm
\centerline{\hh ON CERTAIN LARGE ADDITIVE FUNCTIONS}

\bigskip
\centerline{\cc Aleksandar Ivi\'c}

\bigskip\no

\centerline{\sevenbf Dedicated to the memory of Paul Erd\H os
(1913-1996)}
\bigskip
\bigskip\no
{\ss Abstract}. {\rr Let $P(n)$ denote the largest prime factor of
an integer $n \ge 2$, $P(1) = 1$, and let $$ \b(n) =
\sum_{p|n}p,\;B(n) = \sum_{p^\a||n}\a p,\;B_1(n) = \sum_{p^\a||n}
p^\a $$ denote ``large" additive functions. A survey of results on
these functions is presented, as well as some new results and open
problems.}

\bigskip\no
\centerline{\dd 1. Introduction}
\bigskip\no

Let $P(n)$ denote the largest prime factor of an integer $n \ge 2$,
$P(1) = 1$, and let
$$
\b(n) = \sum_{p|n}p,\;B(n) = \sum_{p^\a||n}\a p,\;B_1(n) = \sum_{p^\a||n}
p^\a\leqno(1.1)
$$
denote ``large" additive functions, in contrast with the well-known
``small" additive functions
$$\omega(n) = \sum_{p|n}1,\quad
\Omega(n) = \sum_{p^\a||n}\a.\leqno(1.2)
$$
As usual $p$ will denote primes,
$p^\a||n$ means that $p^\a$ divides $n$ but $p^{\a+1}$ does not, and
a function $f(n)$ is {\it additive} if $f(mn) = f(m) + f(n)$ whenever
$(m,n)=1$.  From the pioneering works of Alladi
and Erd\H os [1]-[2], P. Erd\H os's perspicacity and insight have been
one of the main driving forces in the research that brought on many
results on summatory functions of large additive functions and $P(n)$.
The functions $\omega(n)$ and $\Omega(n)$ may be successfully
investigated by various analytical methods. In fact, it is Erd\H os
who in two classical works with M. Kac [14], [15] established the
Gaussian distribution law for these functions. A general principle
is that $z^{f(n)}$ is a multiplicative function whenever $f(n)$
is an additive function. Thus from the Euler product representation
$$
\sum_{n=1}^\infty z^{\omega(n)}n^{-s}
= \prod_p\left(1 + {z\over p^s-1}\right)
= \zeta^z(s)G(s,z)\quad ({\rm Re\,}s > 1, z \in \CC),\leqno(1.3)
$$
where the Dirichlet series for $G(s,z)$ converges absolutely for
${\rm Re\,}s > {1\over2}$, one can obtain various results involving
the distribution of values of $\omega(n)$ (and similarly
of  $\Omega(n)$ and
$\Omega(n) - \omega(n)$; see e.g., [22, Chapter 13]). However, from
the analogue of (1.3) for $\b(n)$, namely
$$
\sum_{n=1}^\infty z^{\b(n)}n^{-s}
= \prod_p\left(1 + {z^p\over p^s-1}\right)\qquad
({\rm Re\,}s > 1, |z| \le 1) \leqno(1.4)
$$
one cannot factor out a power of $\zeta(s)$ that will dominate the Euler
product (because of the difficulties inherent in handling the factor
$z^p$), as was the case in (1.3). Therefore (1.4) does not appear
to be very useful in dealing with problems involving $\b(n)$.

\medskip
For this reason other methods of approach seemed more appropriate to
use. They involve a combination of various analytic and elementary methods.
It transpired
that in many problems a decisive r\^ole is played by the function
$$
\psi(x,y) \;=\; \sum_{n\le x,P(n)\le y}1,\leqno(1.5)
$$
the number of integers not exceeding $x$ all of whose prime factors
do not exceed $y$. Results by Hildebrand, Tenenbaum (see [17]--[19] and
[34]) and others brought on great progress.  In more ways than one this
progress on $\psi(x,y)$ is reflected on the results on
large additive functions and the largest prime factor of an integer.

\medskip
The purpose of this paper is to present an overview of some of
the results on large additive functions and the largest prime
factor of an integer. This topic is motivated by the joint
works of P. Erd\H os and the author [6]--[13],
where the majority of published papers deals precisely with
large additive functions and $P(n)$. The span of the research covers
a period of more than fifteen years, and besides P. Erd\H os and
the author involves works of J.-M. De Koninck [3]-[5], C. Pomerance
[15], [26], Smati and Wu [32], [33], Tizou Xuan [35], [36] and others.
As already mentioned, it was  P. Erd\H os who was
the driving force behind this research, always ready to listen to
ideas and problems, and always prepared to pour out new problems
of his, new methods, and new ideas.

\bigskip\no
\centerline{\dd 2. Some results on summatory functions}
\bigskip\no
The first results on the summatory functions of $\b(n)$ and $B(n)$
we obtained by Alladi--Erd\H os [1], [2]. Later research refined
some of their results. Now we know that $$\eqalign{ \sum_{n\le
x}\b(n) &\;=\; \sum_{j=1}^MA_j{x^2\over\log^jx} +
O\left({x^2\over\log^{M+1}x}\right),\cr A_j &\;=\;
(-1)^{j-1}\,{\d^{j-1}\{\zeta(s)s^{-1}\}\over\d s^{j-1}}
\Bigg|_{s=2}, \qquad A_1 = {\pi^2\over12},\cr} \leqno(2.1) $$ for
any fixed integer $M \ge 1$ (see [4], and [5] for the analogues
for large additive functions over primes of positive density). The
asymptotic formula (2.1) holds if $\b(n)$ is replaced by $P(n)$ or
$B(n)$, and it also holds if $\b(n)$ is replaced by $B_1(n)$,
since one has $$ \sum_{n\le x}B_1(n) \;=\;\sum_{n\le x}\b(n) +
O(x^{3/2}). $$ The summatory functions of quotients of large
additive functions were extensively investigated. The work of P.
Erd\H os and the author [8]  contains proofs of $$\eqalign{
 \sum_{2\le
n\le x}{f(n)\over B_1(n)}
 &= x + O\left({x\log\log x\over\log x}
\right) \qquad(f(n) \in \{P(n),\,\b(n),\,B(n)\}),\cr
 \sum_{2\le
n\le x}{B_1(n)\over B(n)} &= Dx + O\left({x\over\log^{1/3}x}
\right) \;(D > 0),\cr
 \sum_{2\le n\le x}{ B_1(n)\over g(n)} &=
e^\gamma x\log\log x + O(x),\cr}
$$
where $\gamma$ is Euler's
constant, and $g(n) \in \{P(n),\b(n)\}$. The ``closeness" of
$\b(n)$ and $B(n)$ is also evident in the asymptotic formula (see
[25]) $$ \sum_{n\le x}\big(B(n)-\b(n)\big) = x\log\log x +
x\sum_{j=0}^M {C_j\over\log^jx}  +
\,O\left({x\over\log^{M+1}x}\right), $$ which is valid for any
fixed integer $M \ge 1$ and suitable constants $C_j$.

\bigskip
\centerline{\dd 3. Local densities of $B(n) - \b(n)$}
\bigskip\no
The ``local density" of a nonnegative, integer-valued
arithmetic function $f(n)$ is the quantity
$$
d_k \;:=\; \lim_{x\to\infty}\,{1\over x}\,\sum_{n\le x, f(n)=k}1
\qquad (k \in \NN\cup\{0\}),
$$
provided that the limit exists. A classical problem of analytic
number theory are the local densities of $\Omega(n) - \omega(n)$,
which is known as ``R\'enyi's problem" (see [22, Chapter 13]).
For a discussion on local densities of a fairly wide class of
arithmetic functions the reader is referred to Ivi\'c--Tenenbaum [27].
Here we shall complement the results on  $B(n)$ and $\b(n)$ by presenting
a new result. This is the following

\medskip\no
THEOREM. {\it For suitable constants $d_k$ we have, uniformly
in $k \ge 0$},
$$
\sum_{n\le x,B(n)-\b(n)=k}1 \;=\; d_kx \;+\; O(x^{1\over2}\log x),\leqno(3.1)
$$
$$\sum_{n\le x,B(n)-\b(n)\ge k}1 \;\ll\; {x\over k}\qquad(k \ge 1)
,\leqno(3.2)
$$
{\it and also for any given} $r > 0$
$$
{\mathop{\sum\nolimits'}_{n\le x}}{1\over(B(n)-\b(n))^r} =
D_rx + O(x^{1-{r\over2}}\log x) + O(x^{1\over2}\log x), \leqno(3.3)
$$
{\it where $\sum\nolimits'$ denotes summation over  numbers which are
not squarefree and}
$$
D_r \;=\; \sum_{k=1}^\infty {d_k\over k^r}.\leqno(3.4)
$$
\smallskip\no {\bf Proof.} The asymptotic formula in (3.1)
follows similarly as the proof
of the author's result [20] (see  also [22, Chapter 13])
for the function  $a(n)$ (the number of non-isomorphic Abelian
groups with $n$ elements). One has only to replace $a(n)$
by $B(n)-\b(n)$, since
$$
B(pm)-\b(pm) = B(m)-\b(m)
$$
if $(p,m) = 1$, and the method of proof of [20] goes through.
If $q$ denotes squarefree and $s$ denotes squarefull numbers, then
since every $n$ can be written uniquely as $n = qs,\, (q,s) = 1$,
the  sum in (3.2) becomes
$$
\sum_{qs\le x,(q,s)=1,B(s)-\b(s)\ge k}1 \le x\sum_{s=1, B(s)-\b(s)\ge k}
^\infty\,{1\over s}.
$$
But by induction on $\omega(s)$ one has
$$
k \;\le\; B(s) - \b(s) \;=\; \sum_{p^\a||s}(\a-1)p \;\le\;
 2\prod_{p^\a||s}p^{\a\over2} \;=\;2\sqrt{s},\leqno(3.5)
$$
hence $s \ge k^2/4$, which because of $\sum_{s\le x}1 \ll \sqrt{x}$ gives
$$
\sum_{s=1, B(s)-\b(s)\ge k}^\infty\,{1\over s} \;\le\;
 \sum_{s\ge k^2/4}
\,{1\over s} \;\ll\; {1\over k}.
$$
To prove (3.3) note that the sum on the left hand-side is, in view
of (3.1) and (3.5),
$$\eqalign{&
\sum_{1\le k\le2\sqrt{x}}{1\over k^r}\left(\sum_{n\le x,B(n)-\b(n)=k}1\right)
= \sum_{1\le k\le2\sqrt{x}}{1\over k^r}\left(d_kx
+ O(x^{1\over2}\log x)\right)
\cr& = \sum_{k=1}^\infty{d_k\over k^r}\,x + O\left(
x\sum_{k>2\sqrt{x}}k^{-1-r}\right)
+ O\left(x^{1\over2}\log x\sum_{k\le2\sqrt{x}}k^{-r}\right)\cr&
= D_rx + O\left(x^{1-{r\over2}}\log x\right) + O(x^{1\over2}\log x),\cr}
$$
with $D_r$ as in (3.4).

\medskip We shall conclude this section by stating some open problems.

\smallskip
{\bf Problem 1}. Which density $d_k$ is the largest one for $k > 1$?
(We have $d_0 = 6/\pi^2$ (= the density of squarefree numbers),
$d_1 = 0$ and $d_k > 0$ for $k > 1$).

\medskip
{\bf Problem 2}. The  proof of (3.2) shows that $d_k \ll 1/k$. Is $d_k$
decreasing for $k \ge k_0$? Can one find an asymptotic formula for
$d_k$?

\medskip
{\bf Problem 3}. Is the density of $n$ for which $\b(n) > \b(n+1)$ (or
$B(n) > B(n+1), B_1(n) > B_1(n+1)$) equal to 1/2? What about the density
of $n$ for which, say, $\b(n) > \b(n+1) > \b(n+2)$?

\smallskip\no These are the analogues of Erd\H os's classical problem
to prove that the density of $n$ for which $P(n) > P(n+1)$ is 1/2.

\medskip
{\bf Problem 4}. For which $n$ is it possible to have $\b(n) =\b(n+1)$
(like $\b(5) = \b(6)$),
$B(n) = B(n+1)\;(B(714) = B(715))$ and $B_1(n) = B_1(n+1)$? It was proved by
Erd\H os--Pomerance [16] that
$$
\sum_{n\le x,B(n)=B(n+1)}1 = O\left({x\over\log x}\right).\leqno(3.6)
$$
One could look either for asymptotic estimates such as (3.6), or
try to give an arithmetic characterization of the numbers in question.
\medskip
{\bf Problem 5}. Can one improve the $O$--term in (3.1) by taking into
account the arithmetic structure of $k$?

\smallskip\no It may be remarked that in the analogous problem
for the local densities of $a(n)$
(the number of non-isomorphic Abelian groups with $n$ elements) this was done
by Kr\"atzel--Wolke [28].

\bigskip
\bigskip
\centerline{\dd 4. Sums of reciprocals}
\bigskip\no
It is a classical result of prime number theory that
$$
\sum_{p\le x}{1\over p} = \log\log x + C + O\left({1\over\log x}\right).
$$
Sums of reciprocals of large additive functions and $P(n)$  are much
more difficult to handle. They were
investigated by De Koninck, Erd\H os, Pomerance, Xuan and the author.
It was proved by  Erd\H os, Pomerance and the author [13] that
$$
\sum_{n\le x}{1\over P(n)} = x\delta(x)\left(1 + O\left(\sqrt{{\log\log x
\over\log x}}\,\right)\right)\leqno(4.1)
$$
with
$$ \delta(x) \;:=\; \int_2^x \rho\left(
{\log x\over\log t}\right){\d t\over t^2},\leqno(4.2)
$$
where the Dickman--de Bruijn function $\rho(u)$
is the continuous solution
to the differential delay equation
$$u\rho'(u) = -\rho(u-1), \;\rho(u) = 1\;
{\rm for}\; 0 \le u \le 1,\; \rho(u) = 0\;{\rm for}\; u < 0\;.$$
It is known (see [34]) that ($\log_kx = \log(\log_{k-1}x$))
$$
\rho(u) = \exp\Biggl\{-u\Biggl(\log u + \log_2u - 1
+ {\log_2u-1\over\log u} + O\left(\left({\log_2u\over\log u}\right)^2
\right)\Biggr)\Biggr\}.
$$
It was also proved in [13] that
$$
\sum_{n \le x}P(n)^{-\omega(n)} = \exp\left\{(4+o(1)){\sqrt{\log x}
\over\log\log x}\right\},\leqno(4.3)
$$
$$
\sum_{n \le x}P(n)^{-\Omega(n)} = \log\log x + D +
O\left({1\over\log x}\right)\quad (D > 0),\leqno(4.4)
$$
with effectively computable $D$, showing the difference in behaviour between
$\omega(n)$ and $\Omega(n)$. Of these two formulas it is (4.3) that
is deeper than (4.4).

\medskip
{\bf Problem 6.} What is the shape of the above asymptotic formulas
if we replace $P(n)$ by $\b(n)$ and $B(n)$?

\medskip
It was shown by Pomerance and the author [26] that one has
asymptotically $$ \delta(x) = \exp\left\{-(2\log x\log_2x)^{1/2}
\left(1 + g_0(x) +
O\left({\log_3^3x\over\log_2^3x}\right)\right)\right\}, $$ where
$$ g_r(x) = {\log_3x + \log(1+r) -2 -\log2\over2\log_2x}\left(1 +
{2\over\log_2x}\right) $$ $$- {(\log_3x + \log(1+r) -
\log2)^2\over8\log_2^2x}, $$ and the expression for $\delta(x)$
was sharpened by the author [24]. Already in 1977 Erd\H os told
the author that the function $\delta(x)$ is slowly varying in the
sense of J. Karamata (see [29], [31]), namely that for any $C > 0$
one has $$ \lim_{x\to\infty}\;{\delta(Cx)\over\delta(x)} \;=\;
1,\leqno(4.5) $$ but it is only in 1986 in that (4.5) was
established in [13], by the use of (4.2) and properties of the
function $\rho(u)$. The asymptotic formula (4.1) remains valid if
$P(n)$ is replaced by $\b(n)$ or $B(n)$, and the asymptotic
formula for the summatory function of ${B(n)\over\b(n)}-1$ is of
the same shape as the right-hand side of (4.1). Furthermore we
have $$\eqalign{& \sum_{2\le n\le x}\left({1\over\b(n)} - {1\over
B(n)}\right) \cr& =  x\exp\left\{-2(\log x\log_2x)^{1/2}\left(1 +
g_1(x) +
O\left({\log_3^3x\over\log_2^3x}\right)\right)\right\}.\cr} $$
Based on his joint work with   Erd\H os and  Pomerance [13], the
author [23] sharpened some of the  formulas and obtained e.g.
$$\eqalign{
 \sum_{n\le x}{\Omega(n)-\omega(n)\over P(n)} &=
\left\{\sum_p{1\over p^2-p}+ O\left(\left({\log_2x\over\log
x}\right) ^{1/2}\right)\right\}\sum_{n\le x}{1\over P(n)}, \cr
\sum_{n\le x}{\omega(n)\over P(n)} &=
 \left\{\left({2\log x\over\log_2x}\right)^{1/2}
\left(1 +  O\left({\log_3x\over\log_2 x}\right)
\right)\right\}\sum_{n\le x}{1\over P(n)},\cr}
 $$
  and this remains
valid if $\omega(n)$ is replaced by $\Omega(n)$, $$ \sum_{n\le
x}{\mu^2(n)\over P(n)} =
\left\{{6\over\pi^2}+O\left(\left({\log_2x\over\log
x}\right)^{1/2} \right)\right\}\sum_{n\le x}{1\over P(n)}. $$ The
relevant contribution to the last three sums comes from $n$ for
which $$ L(-2,x) \;\le\; P(n) \;\le\;  L(2,x), $$ where $$ L(c,x)
\;:=\; \exp\left\{({\textstyle{1\over2}}\log x\log_2x)^{1/2}
\left(1 + c{\log_3x\over\log_2x}\right)\right\}. $$ During many
years of collaboration on the problems discussed in this section
Erd\H os was fond of saying ``take $\log P(n) \asymp \sqrt{\log
n}\;$", but the above discussion shows that in this, as on
countless other occasions, he was right.

\bigskip
\centerline{\dd 5. Sums in residue classes}
\bigskip\no
During the Conference on Analytic Number Theory in June 1993 in
Lillaf\"ured, P. Erd\H os asked the author
to evaluate asymptotically $S_0(x)$,
where for fixed $r \ge 0$ and fixed integers $1 \le \ell \le k, (\ell,k)
= 1$ one defines
$$
S_r(x) \;=\; \sum_{n\le x,P(n)\equiv\ell({\rm mod})k}{1\over P^r(n)}.
$$
In the author's work [24] it is shown that that
$$
S_0(x) = {x\over\varphi(k)} + O\left(x\exp(-(\log x)^{3/8-\e})\right)
$$
for any given $\e > 0$, and for $r > 0$ and any fixed integer $J \ge 0$
$$
S_r(x) \;=\; {x\over\varphi(k)}\int_2^x\rho\left({\log x\over\log t}\right)
\Biggl\{\sum_{j=0}^J{Q_{j,r}(\log t)\over\log^jx} \,+
O\left(\left({\log t\over\log x}\right)^{J+1}\right)\Biggr\}
{\d t\over t^{r+1}}\leqno(5.1)
$$
for suitable polynomials $Q_{j,r}(x)$ of degree $j$ in $x$ whose
coefficients depend on $r$. In particular
$$
Q_{0,r}(x) = r,\quad Q_{1,r}(x) = (r - r\gamma)(rx - 1).
$$
Let
$$
T_r(x) := \sum_{n\le x,P(n)\equiv \ell({\rm mod}\,k),P^2(n)|n}
{1\over P^r(n)}\quad(r \ge -1),
$$
$r \in \RR$ fixed, $1 \le \ell\le k,\,(\ell,k) = 1$  fixed. For
any given $\varepsilon > 0$
$$ \eqalign{
T_{-1}(x) &\;=\; {Cx\over\varphi(k)} + O\left\{x\exp\left(-\log^{{3\over8}
-\varepsilon}x\right)\right\}, \cr
 C &\;=\; \int_0^\infty {\rho(v)\over v+2}\d v\quad(< 1). \cr}
$$
For $r > -1, \, r\in\RR$ fixed and $J\in\NN$ fixed
$$
T_r(x) =  {x\over\varphi(k)}
\int_2^x\rho\left({\log x\over\log t}\right)
\left(\sum_{j=0}^J{R_{j+1,r}(\log t)\over \log^jx} +
O\left({\log^{J+2}t\over\log^{J+1}x}\right)\right)
{\d t\over t^{r+2}}\leqno(5.2)
$$
for suitable polynomials $R_{j,r}(x)\;(j \in \NN)$ of degree $j$ in
$x$ whose coefficients depend on $r$. In particular,
$$
R_{1,r}(x) \;=\; (r+1)^2x - r -1.
$$
The formulas (5.1) and (5.2) sharpen the results of [26] (in the
case when $k = 1$), where one had
$$
S_r(x) = x\exp\Biggl\{-(2r\log x\log_2x)^{1/2}\Biggl(1 + g_{r-1}(x) +
 O\left(\left({\log_3x\over\log_2x}\right)^3\right)
\Biggr)\Biggr\},
$$
when $r > 0$, and
$$
T_r(x) = x\exp\Biggl\{-((2r+2)\log x\log_2x)^{1/2}
\Biggl(1 + g_{r}(x) +
O\left(\left({\log_3x\over\log_2x}\right)^3\right)
\Biggr)\Biggr\},
$$
when $r > -1$.
The proofs given in [24] use the sharp approximation of E. Saias [30], namely
$$
\psi(x,y) =  \Lambda(x,y)\left\{1 + O\left(\exp\left(-\log^{{3\over5}-
\varepsilon}y\right) \right) \right\},\leqno(5.3)
$$
where
$$e^{(\log\log x)^{5/3+\varepsilon}} \le y \le x,\; x \ge
x_0(\varepsilon),\leqno(5.4)
$$
for $ \; x,y \ge 1,\,x\not\in\NN $
$$
\Lambda(x,y) := x\int_{1-0}^\infty\rho\left({\log x-
\log t\over\log y}\right)\d\left({[t]\over t}\right),\leqno(5.5)
$$
and for $x\in\NN$ we define $ \Lambda(x,y) = \Lambda(x+0,y)$.
Note that in the range (5.4) A. Hildebrand [17] had
$$
\psi(x,y) = x\rho(u)\left(1 + O\left({\log(u+2)\over\log y}\right)\right),
\quad u = {\log x\over \log y}.\leqno(5.6)
$$
Although (5.3) has the sharper error term than (5.6),
the function on the right-hand side of (5.5) is
discontinuous ($[t]$, the greatest integer part of $t$ has jumps
when $t\in\NN$), and there are technical
difficulties in applying this formula.

\smallskip
As one of the corollaries of the above results we single out
the following formula:
$$
\sum_{n\le x,P^2(n)|n}1 = \Biggl\{{\log x\over2}\Biggl(\log_2x +
\log_3x - \log2 \,+
$$
$$
+\; {\log_3x-\log2\over\log_2x}+ O\left({\log_3^2x\over\log_2^2x}
\right)\Biggr)\Biggr\}^{1/2}\,\sum_{n\le x}{1\over P(n)}.
$$
The error term, like in most previous results, could be further
sharpened at the cost of more technical elaboration.

\vfill\eject \topglue2cm
 \centerline{\cc References}

\bigskip\sevenrm
\item{[1]} K. Alladi and P. Erd\H os, On an additive arithmetic function,
Pacific. J. Math. {\sevenbf71}(1977), 275-294.
\item{[2]} K. Alladi and P. Erd\H os, On the asymptotic behavior of
large prime factors of integers, Pacific J. Math. {\sevenbf82}
(1979), 295-315.

\item{[3]} J.-M. De Koninck and A. Ivi\'c, Topics in arithmetical
    functions, Mathematics Studies {\sevenbf43}, North-Holland,
    Amsterdam 1980.

\item{[4]} J.-M. De Koninck and A. Ivi\'c,
 The distribution of the average prime
divisor of an integer, Archiv  Math.
{\sevenbf43}(1984), 37-43.

\item{[5]} J.-M. De Koninck and A. Ivi\'c,
On some asymptotic formulas related
    to large additive functions over primes of positive density,
    Mathematica Balkanica  {\sevenbf10} (1996), 279-300.

\item{[6]} J.-M. De Koninck, P. Erd\H os and A. Ivi\'c,
    Reciprocals of certain large additive
    functions, Canadian Math. Bulletin {\sevenbf24} (1981), 225-231.

\item{[7]}   P. Erd\H os, S.W. Graham,  A. Ivi\'c and C. Pomerance,
       On the number of divisors of n!,
Analytic Number Theory: Proceedings of a
       Conference in Honor of Heini Halberstam (Urbana, May 1995) Volume 1,
(eds. B.C. Berndt et al.), Birkh\"auser, Boston etc., 1996, 337-355.

\item{[8]} P. Erd\H os and A. Ivi\'c,
Estimates for sums involving the largest prime     factor of an
    integer and certain related additive functions, Studia Scientiarum
    Math. Hungarica {\sevenbf15} (1980), 183-199.

\item{[9]}  P. Erd\H os and A. Ivi\'c,
 On sums involving reciprocals of certain arithmetical
functions,
Publications Inst. Math. (Belgrade) {\sevenbf 32(46)} (1982), 49-56.

\item{[10]}  P. Erd\H os and A. Ivi\'c,
 The distribution of certain arithmetical functions
    at consecutive integers,
Proceedings Budapest Conference in Number Theory July 1987,
Coll. Math. Soc. J. Bolyai {\sevenbf 51}, North-Holland, Amsterdam 1989,
    45-91.

\item{[11]}  P. Erd\H os and A. Ivi\'c,
On the iterates of the enumerating function of
finite Abelian groups,  Bulletin XCIC Acad.
Serbe 1989 Sciences Math\'ematiques    No {\sevenbf 17}, 13-22.

\item{[12]}  P. Erd\H os and A. Ivi\'c,
 The distribution of small and large additive functions II,
 Proceedings of the Amalfi Conference on Analytic Number Theory
    (Amalfi, Sep. 1989), Universit\`a di Salerno, Salerno 1992, 83-93.

\item{[13]} P. Erd\H os, A. Ivi\'c and C. Pomerance,
On sums involving reciprocals of the
largest prime factor of an integer,  Glasnik Matemati\v cki
{\sevenbf21(41)} (1986), 283-300.

\item{[14]} P. Erd\H os and M. Kac,
On the Gaussian law of errors in the theory
of additive functions, Proc. Nat. Acad. Sc. U.S.A.
{\sevenbf25}(1939), 206-207.

\item{[15]} P. Erd\H os and M. Kac, The Gaussian law of errors  in the theory of
additive number-theoretic functions, Amer. Journal Math.
{\sevenbf62}(1940), 738-742.

\item{[16]} P. Erd\H os and C. Pomerance, On the largest prime factors
of {\seveni n} and {\seveni n}+1, Aequationes Math. {\sevenbf17}(1978), 311-321.

\item{[17]} A. Hildebrand, On the number of positive integers $\le$
{\seveni x} and free of prime factors $>$ {\seveni y}, J. Number Theory
{\sevenbf22}(1986), 289-307.

\item{[18]} A. Hildebrand and G. Tenenbaum, On integers free of large
prime factors, Trans. Amer. Math. Soc. {\sevenbf 296}(1986), 265-290.

\item{[19]} A. Hildebrand and G. Tenenbaum, Integers wuthout large prime
factors, J. Th\'eorie des Nombres Bordeaux {\sevenbf5}(1993), 411-484.

\item{[20]} A. Ivi\'c, The distribution of values of
the enumerating function of non-isomorphic
   abelian groups of finite order,  Archiv  Math. (Basel-Stuttgart)
   {\sevenbf30}    (1978), 374-379.

\item{[21]} A. Ivi\'c, Sum of reciprocals of the largest
prime factor of an integer, Archiv
    Math. {\sevenbf36} (1981), 57-61.

\item{[22]} A. Ivi\'c, The Riemann zeta-function, John Wiley \& Sons,
New York, 1985.

\item{[23]} A. Ivi\'c, On some estimates involving the
number of prime divisors of an integer,
    Acta Arithmetica {\sevenbf49} (1987), 21-32.

\item{[24]} A. Ivi\'c, On sums involving reciprocals of
the largest prime factor of an integer II,
    Acta Arithmetica {\sevenbf71} (1995), 241-245.

\item{[25]} A. Ivi\'c, On large additive functions over
primes of positive density,
    Mathematica Balkanica {\sevenbf10} (1996), 103-120.

\item{[26]} A. Ivi\'c and C. Pomerance,
 Estimates for certain sums involving the
    largest prime factor of an integer, Proceedings Budapest Conference
    in Number Theory July 1981, Coll. Math. Soc. J. Bolyai {\sevenbf34},
    North-Holland, Amsterdam 1984, 769-789.

\item{[27]} A. Ivi\'c and G. Tenenbaum, Local densities over integers free
of large prime factors, Quart. J. Math. (Oxford) (2) {\sevenbf 37}(1986),
401-417.

\item{[28]} E. Kr\"atzel and D. Wolke,
\"Uber die Anzahl der Abelschen Gruppen
gegebener Ordnung, Analysis {\sevenbf14}(1994), 257-266.

\item{[29]} J. Karamata, Sur un mode de croissance r\'eguli\`ere des
fonctions, Mathematica (Cluj) {\sevenbf4}(1930), 38-53.

\item{[30]} E. Saias, Sur le nombre des entiers sans grand facteur
premier, J. Number Theory {\sevenbf 32}(1989), 78-99.

\item{[31]} E. Seneta, Regularly varying functions, LNM 508, Springer
Verlag, Berlin--Heidelberg--New York, 1976.

\item{[32]} A. Smati, Sur l'it\'eration du nombre de diviseurs des
entiers sans grand facteur premier,  J. Number Theory
{\sevenbf 57}(1996), 66-89.

\item{[33]} A. Smati and J. Wu, Distribution of Euler's function
over integers free of large prime factors, Acta Arith.
{\sevenbf 72}(1996), 139-155.

\item{[34]} G. Tenenbaum, Introduction \`a la th\'eorie analytique
et probabiliste des nombres, Soci\'et\'e Math. de France, Paris, 1995.

\item{[35]} T.Z. Xuan, On sums involving reciprocals of certain large
additive functions, Publs. Inst. Math. (Beograd) {\sevenbf45(59)}, 41-55
and II, ibid. {\sevenbf46(60)}(1989), 25-32.

\item{[36]} T.Z. Xuan, On a result of Erd\H os and Ivi\'c, Archiv Math.
{\sevenbf62}(1994), 143-154.

\bigskip\bigskip
\no Aleksandar Ivi\'c

\no Katedra Matematike RGF-a

\no Universiteta u Beogradu

\no Dju\v sina 7, 11000 Beograd

\no Serbia (Yugoslavia)

\no e-mail: \sevenbf aivic@matf.bg.ac.yu, aivic@rgf.bg.ac.yu

\bye